\input amstex.tex
\documentstyle{amsppt}
\input xy
\xyoption{all}
\CompileMatrices
\xyReloadDrivers

\def\t#1 {\text{\rm #1}}
\def\b#1 {\text{\rm\bf #1}}


\topmatter
\title A controlled-topology proof of the product structure theorem
\endtitle
\author Frank Quinn\endauthor
\address Mathematics, Virginia Tech, Blacksburg VA 24061-0123\endaddress
\email quinn\@math.vt.edu\endemail \date September 2006\enddate
 \subjclass 57Q25\endsubjclass

\abstract 
The controlled end and h-cobrodism theorems (Ends of maps I, 1979) are used to give quick proofs of the Top/PL and PL/DIFF product structure theorems. 
\endabstract
\endtopmatter

\head 1. Introduction\endhead
The ``Hauptvermutung'' expressed the hope that a topological manifold might have a unique PL structure, and perhaps analogously a PL manifold might have a unique differentiable structure. This is not true so the real theory breaks into two pieces: a way to distinguish structures; and the proof that this almost always gives the full picture. 

Milnor's microbundles \cite{M} are used to distinguish structures. This is a relatively formal theory. Let $\Cal M\subset \Cal N$ denote two of the manifold classes $\t DIFF \subset\t PL  \subset\t TOP $. An $\Cal N$ manifold $N$ has a stable normal (or equivalently, tangent) $\Cal N$ microbundle, and an $\Cal M$ structure on $N$ provides a refinement to an $\Cal M$ microbundle. The theory is set up so that this automatically gives a bijective correspondence between stable refinements of the bundle and stable structures, i.e\. $\Cal M$ structures on the $\Cal N$ manifold $N\times R^k$, for large $k$. 

Microbundles can be described using classifying spaces. This shows, for instance, that if $N$ is a topological manifold then $M\times R^k$ has a smooth structure (for large $k$) exactly when the stable tangent (or normal) microbundle $N\to B_{\t TOP }$ has a (homotopy) lift to $B_{\t DIFF }$, and stable equivalence classes of such structures correspond to homotopy classes of lifts. The homotopy type of the classifying spaces can be completely described modulo the usual mystery of stable homotopy of spheres.

The geometrically difficult step is the ``product structure theorem'' that describes when structures on $N\times R^k$ correspond to structures on $N$. The Cairns--Hirsh theorem asserts this is always true for Diff $\subset$ PL. The Kirby--Siebenmann theorem is that the PL $\subset$ TOP version is true for manifolds of dimension greater than 5. Donaldson showed this is false in dimension 4. 

The breakthrough for the topological product structure theorem was due to Kirby \cite{K} and extended by Siebenmann \cite{KS1} in 1969. The approach was long and elaborate and  a full account did not appear until 1977 \cite{KS2}. A more direct proof using controlled topology became available only two years later \cite{Q1} and is described here. This approach has been effective in further investigations. It was used in  \cite{Q2} to describe obstructions to existence and uniqueness of PL structures on general locally triangulable spaces, including a  counterexample (in fact the quotient of a locally linear finite group action on $D^7$) so that the complement of the singular set does not have the homotopy type of a finite complex.  It was used in \cite{Q3},  see also \cite{FQ \S8.1}, to show that some fragments of the theory (the ``annulus conjecture'' and microbundle stability) remain valid in dimension 4. 

The controlled approach is more direct in the sense that the product structure theorem is an easy deduction from other results, but these other results are far from elementary. The hard work is  shifted to the proof of a vanishing theorem for obstruction groups. The real payoff is not overall ease but a wealth of other applications.  

\head 2. Product structure theorem\endhead
Let $\Cal M\subset \Cal N$ denote two of the manifold classes $\t DIFF \subset\t PL  \subset\t TOP $. The primitive definition of a ``structure'' in the PL or DIFF cases is in terms of systems of coordinate charts on a fixed space. However in making comparisons it is easier to think of the structures as being on different spaces. 
\proclaim{2.1 Definition} \roster \item An $\Cal M$ structure on an $\Cal N$ manifold $N$ is an $\Cal M$ manifold $M$ and an  $\Cal N$  isomorphism $M\to N$.
\item A concordance (between two structures) is a structure on $N\times I$ (that restricts to the given structures over the ends).
\item a structure is ``rel boundary'' if $\partial N$ already has an $\Cal M$ refinement in the primitive chart sense  and $M\to N$ is an $\Cal M$ isomorphism.
\endroster
\endproclaim
Mapping cylinders of isotopies provide special examples of concordances. Suppose $\tau\colon M_0\to M_1$ is an $\Cal M$ isomorphism, $\theta_i\colon M_i\to N$ are $\Cal N$ structures, and $h\colon \theta_1\tau \sim \theta_0$ is an $\Cal N$ isotopy. Let $\t Cyl (\tau)= M_0\times[0,1]\cup M_1/\simeq$ denote the
 mapping cylinder, with equivalence relation defined by $(x,1)\simeq\tau(x)$. Then the map $(x,t)\mapsto (h(x,t),t)$ is an $\Cal N$ isomorphism and therefore a concordance.

\proclaim{2.2 Product Structure Theorem} Suppose
 $N$ is an $\Cal N$ manifold, $\partial N$ has an $\Cal M$ structure, and the dimension of $N$ is $\geq 5$. Then:
\par\noindent{\bf Stability:}\quad If\/ $\theta\colon M\to N\times R^k$ is an $\Cal M$ structure rel boundary then there is an $\Cal M$ structure $\hat \theta\colon \hat M\to N$  and an $\Cal M$ concordance rel boundary between $\theta$ and\/ $\hat\theta\times \t id \colon\hat M\times R^k\to N\times R^k$.
 
 \par\noindent{\bf Isotopy:}\quad Suppose $\theta\colon M\to N\times I$ is a concordance  rel boundary between  $\Cal M$ structures $\partial_i M_i\to N\times\{i\}$, $i=0,1$.  Then  there is an
 $\Cal M$ isomorphism rel boundary from a mapping cylinder  $\beta\colon \t cyl (\tau)\to M$ where $\tau\colon M_0\to M_1$ is an $\Cal M$ isomorphism rel boundary; and
 $\Cal N$ isotopies $h\colon(\partial_1\theta)\tau\sim \partial_0\theta$ and  $\theta\beta\sim \t cyl (h)$.
  \endproclaim

\subhead 2.3 Notes\endsubhead
\roster
\item The Isotopy statement asserts that concordant structures are isomorphic, and the isomorphisms are well-defined up to isotopy. The isomorphisms and isotopies can be chosen smaller than any given $\delta>0$, and this will play an important role. The isomorphisms and isotopies can be thought of as particularly nice concordances; see the comment after 2.1.
\item As mentioned in \S1 this statement remains true without dimension restriction for $\t DIFF \subset \t PL $. Low dimensions can be handled by direct arguments showing $\t DIFF = \t PL $ for manifolds of dimensions $\leq 5$. Parts of this argument are quite long \cite{C, H}. 
\item In dimension 4 the $\t PL \subset \t TOP $ statement is false. However it is possible to get $M_{\infty}$ with a PL structure on the complement of a ``singular set'' with topological dimension 1 and PL dimension  2, and get partial versions from this. See  \cite{Q3} and \cite{FQ \S8.1}. 
\item In dimension 3 topological manifolds all have PL and smooth structures and these are unique up to isotopy. However the product structure theorem as stated is false. First, $\pi_3(TOP/PL)=Z/2$ so the product structure theorem predicts an exotic PL structure on $S^3$. Worse, there are exotic PL (or smooth) structures on $R^4$ that are not products of $R$ with structures on $R^3$.
\item A non rel--boundary version follows from two applications of 2.2. Also, if $U\subset N$ is an open set with an $\Cal M$ structure then a version rel a slightly smaller closed set is obtained by taking a closed codimension--0 $\Cal M$ submanifold of $U$ and then applying 2.2 to the complement of the interior.
\endroster
\head 3. The proof\endhead

 We   refer to \cite{Q1} for statements and explanations of the controlled end and h-cobordism theorems. Quoting statements here would more than double the length of the proof without adding value. Some reductions:
\roster 
\item  For the most part we omit mention of $\partial N$. This simplifies the discussion and putting it  back in is completely routine. 
\item We prove the stability result for $k=1$. The statement is designed so that the general case follows from this by iteration. \endroster

\subhead 3.1 Stability; the structure on $N$\endsubhead
Suppose $\theta\colon M \to N\times R$ is an $\Cal M$ structure. The first step is to take an $\Cal M$ completion over the $+\infty$ end of $N\times R$. A {\it completion\/} in this situation  is an $\Cal M$ manifold with a proper map $\bar\theta\colon\bar M\to N\times(-\infty,\infty]$ so that 
\roster\item $\bar\theta^{-1}(N\times R)=M$, and the restriction of $\bar\theta$ to this is $\theta$; and
\item  $\bar\theta^{-1}(N\times \{\infty\})$ is a subset of $\partial \bar M$, and is denoted $\partial_{\infty}\bar M$. 
\endroster
The hypotheses for finding a completion are that the end is tame and has trivial local fundamental groups over $N$. Since $M\to N\times R$ is cellular and the end of $N\times R$ obviously has these properties, so does the end of $M$. The construction uses handlebody theory in $\Cal M=\t DIFF $ or PL so does not depend on knowing anything about topological manifolds. According to \cite{Q1, Theorem 1.4} an $\Cal M$ completion exists. 

The boundary map $\partial_{\infty}\bar\theta\colon \partial_{\infty}\bar M\to N\times\{\infty\}$ is a map from an $\Cal M$ manifold to $N$.  Since $\theta$ is a homeomorphism this is automatically cellular (point inverses are contractible in arbitrarily small neighborhoods). Modifications needed to get  $\Cal N$ isomorphism divide into  TOP and PL cases. 

First suppose $\Cal N=\t TOP $. A cellular map can be arbitrarily closely approximated by a homeomomorphism. More precisely if $f_0$ is a cellular map of manifolds then there is a continuous 1--parameter family $f_t$ ending in $f_0$ with $f_t$ a homeomorphism for $t>0$. Then $f_1$ in any such 1--parameter family is the structure needed for the theorem.  The homeomorphism approximation of cellular maps when $M$ and $N$ are manifolds of dimensions $\geq5$ was originally due to Siebenbmann \cite{S}. A much more sophisticated theorem that only requires $N$ to be an ANR homology manifold satisfying the disjoint 2--disk property is due to Edwards \cite{D}. A version that only requires the map to be ``approximately'' cellular was given by Chapman and Ferry \cite{CF}. None of these depend on the product structure theorem.

Now suppose $\Cal N=\t PL $, so $\bar M$ is a DIFF manifold. 
The inverse image  $\bar\theta^{-1}(N\times [0,\infty]\,)$ is a PL submanifold of $\hat M$. The restriction of $\hat \theta$ gives a cellular map to $N\times [0,\infty]$, so the projection to $N$ is an $\epsilon$ h-cobordism over $N$ for every $\epsilon>0$. Local fundamental groups are trivial so the controlled h-cobordism theorem \cite{Q1, Theorem 2.7} implies there is a PL trivialization. This gives a PL isomorphism  $\tau\colon \bar\theta^{-1}(N\times [0,\infty]\,)\to N\times [0,\infty]$ that is equal to $\theta$ on $\theta^{-1}(N\times \{0\})$. Denote the restriction to $\partial_{\infty}\bar M$ by $\partial_{\infty}\tau$, then this is a PL isomorphism to $N$, and is the structure needed for the $\Cal N=\t PL $ case. 

\subhead 3.2 Stability; the concordance\endsubhead
The concordances are elementary constructions from the data used to find the structures, essentially by reparameterizing products. 

First suppose $\Cal N=\t TOP $. Reparameterize $(-\infty,\infty]$ as $(-\infty, 0]$ and take the union of $\bar \theta$ and the 1--parameter family $f_t$ of homeomorphisms converging to $\partial_{\infty}\bar\theta$ to get a map 
$$\bar\theta\cup f_*\colon \bar M\cup_{\partial_{\infty}\bar M}\partial_{\infty}\bar M\times [0,1] \to N\times (-\infty, 1].$$
This map is a homeomorphism except over $N\times\{0\}$ where it is cellular, and the domain is an $\Cal M $ manifold. Multiply this by the identity on $[0,1)$ to get a map to $N\times(-\infty,1]\times[0,1)$. Note that over $N\times\{1\}\times[0,1)$ this is the product of the structure $f_1$ with the identity, and the map is a homomorphism except over $N\times \{0\}\times[0,1)$. Delete $N\times[0,1]\times\{0\}$ and its preimage to get an $\Cal M$ manifold and map $\hat M\to N\times\bigl( (-\infty,1]\times[0,1)-[0,1]\times\{0\}\bigr)$. Next we reparameterize the second factor in the range as $(-\infty,\infty)\times[0,1]$ by a homeomorphism that:
\roster\item takes $(-\infty,0)\times\{0\}$ to $(-\infty,\infty)\times \{0\}$ by the inverse of the reparameterization at the beginning of the paragraph; 
\item takes $\{1\}\times (0,1)$ to $(-\infty,\infty)\times \{1\}$; and
\item takes $ \{0\}\times[0,1)$ to $(-\infty,\infty)\times \{1/2\}$.
\endroster
This gives $\hat M\to N\times (-\infty,\infty)\times [0,1]$ that is a homeomorphism except over $N\times (-\infty,\infty)\times \{1/2\}$ where it is cellular. For the final step  approximate this by a homeomorphism unchanged near $N\times (-\infty,\infty)\times \{0,1\}$. This gives the required concordance from $\theta$ to the stabilization of the structure on $N$. 

Now suppose $\Cal N=\t PL $. Recall that the controlled h-cobodism theorem provided a PL isomorphism $\tau$ which extends by $\theta$ on $\theta^{-1}(N\times (-\infty,0]\,)$ to give a PL isomorphism $\tau\colon \bar M\to N\times (-\infty,\infty]$. The plan is to use this to define a map 
$$\bar \tau\colon \bar M\times[0,1)\to N\times(-\infty,\infty]\times [0,1)$$
satisfying:
\roster\item it is $\tau\times \t id $ on $\bar M\times(1/2,1)$;
\item it is $\bar \theta$ on $\bar M\times \{0\}$; and
\item it is a PL isomorphism on $\bar M\times[0,1)-\partial_{\infty}\bar M\times\{0,\}$.
\endroster
The desired concordance is obtained from such a $\bar\tau$ by composing with a reparameterization  $R\times [0,1)\cup \{\infty\}\times (1/2,1)\simeq R\times [0,1]$ similar to the topological case. 

We now define $\bar \tau$. Let $T_x$ denote translation by $x$ on $R$ (i.e\. $T_x(y)=x+y$). Then  for $(x, t)\in \bar M\times[0,1)$
$$\bar \tau(x,t)=\cases (\bar\theta T_{\alpha(t)}\bar\theta^{-1}\hat\theta T_{-\alpha(t)}(x),t)\text{ if }t> 0\\
(\bar\theta(x),0)\text{ if }t=0\endcases$$
where $\alpha\colon (0,1)\to [0,\infty)$ is 0 for $t\geq 1/2$ and increases to $\infty$ as $t\to 0$. For instance $\alpha(t)=\t Max \{\frac1t -2,0\}$. 

This formula is not quite PL because division and translation in $R$ are not PL. It is piecewise smooth, however, so it can be arbitrarily closely approximated by a PL map that is a PL isomorphism except over $N\times(\infty,0)$. 

\subhead 3.3 Isotopy; the isomorphisms\endsubhead
We are given an $\Cal M$ structure $\theta\colon M\to N\times [0,1]$. This is an $\epsilon$ h--cobordism over $N$ for any $\epsilon>0$ and local fundamental groups are trivial so according to \cite{Q1} there is a $\delta$ product structure any $\delta>0$. This product structure is an $\Cal M$ isomorphism $\hat\beta\colon \partial_0M\times [0,1]\to M$ that is the identity on $\partial_0M\times\{0\}$. Restriction to the $\{1\}$ end gives an isomorphism $\tau\colon \partial_0M\to \partial_1M$. There is a canonical isomorphism $\partial_0M\times I\to \t cyl (\tau)$, and $\hat\beta$ factors through this to give the isomorphism $\beta\colon\t cyl (\tau)\to M$ that is the identity on both ends.
\subhead 3.4 Isotopy; the isotopies\endsubhead
Consider the composition $\theta\hat\beta\colon \partial_0M\times I\to N\times I$, where $\hat\beta$ is the $\delta$ product structure from 3.3 just above. We show this is isotopic rel $\partial_0M\times\{0\}$ to an isotopy (i.e\. level-preserving in the $I$ coordinate). This isotopy can be reinterpreted as the two isotopies in the statement of the theorem.

The proof now splits depending on whether $\Cal N$ is PL  or TOP, and we begin with $\Cal N=\t PL $. $\partial_0\theta\colon \partial_0M\to N$ is a PL isomorphism. 
Compose $\theta\hat\beta$ with $(\partial_0\theta)^{-1}\times\t id $ to get a PL isomorphism $N\times I\to N\times I$ that is the identity on $N\times\{0\}$. This is a pseudoisotopy of $N$.  It is sufficient to find an $\Cal N$ isotopy rel $N\times\{0\}$ to the identity. This, however is an immediate application of the controlled pseudoisotopy theorem \cite{Q4, Corollary 1.2}: if $\epsilon>0$ then for sufficiently small $\delta$ (depending on $N$) a PL pseudoisotopy  with trivial local fundamental groups and radius $\leq\delta$ over $N$ is $\epsilon$ PL isotopic to the identity. This completes the proof in this case.

Now suppose $\Cal N=\t TOP $.  As in the PL case compose $\theta\hat\beta$ with $(\partial_0\theta)^{-1}\times\t id $ to get a pseudoisotopy: a homeomorphism $N\times I\to N\times I$ that is the identity on $N\times\{0\}$. Again we want an isotopy to the identity. However according to Edwards--Kirby \cite{EK} the homeomorphism group of $N\times I$ is locally contractible. Thus if $\epsilon>0$ is given there is $\delta>0$ so that a pseudoisotopy within $\delta$ of the identity is isotopic to the identity. Since the argument above gives $\delta$ pseudoisotopies for arbitrarily small $\delta$ the argument is essentially complete. 

There is one last detail. The Edwards--Kirby result requires the pseudoisotopy to be within $\delta$ of the identity when distances are measured in $N\times I$, while the controlled s-cobordism theorem used to produce the pseudoisotopy only gives control in the $N$ coordinate. The Edwards--Kirby proof actually works with this weaker control, but a trick  enables use of the statement without opening up the proof. Recall that the pseudoisotopy is the identity on $N\times\{0\}$. This can be used to find an isotopy that compresses the nontrivial part into a $\delta$ neighborhood of $N\times \{1\}$. The resulting pseudoisotopy is the identity off this neighborhood and the size is unchanged in the $N$ coordinate so it is within $\delta$ of the identity as measured in $N\times I$. This competes the proof in the topological case.

\refstyle{A}

\bye